\documentclass{elsart3-1}

\usepackage{amssymb}

\usepackage[english,francais]{babel}

\newtheorem{theorem}{Theorem}[section]

\newtheorem{lemma}[theorem]{Lemma}
\newtheorem{proposition}[theorem]{Proposition}
\newtheorem{corollary}[theorem]{Corollary}

\newtheorem{conjecture}[theorem]{Conjecture}
\newtheorem{definition}[theorem]{Definition}
\newtheorem{remark}{\it Remark\/}

\renewcommand{\theequation}{\arabic{equation}}
\def\og{\leavevmode\raise.3ex\hbox{$\scriptscriptstyle\langle\!\langle$~}}
\def\fg{\leavevmode\raise.3ex\hbox{~$\!\scriptscriptstyle\,\rangle\!\rangle$}}
\def\EMdash{\leavevmode\hbox to 7.5mm{\vrule height .63ex depth -.59ex width 5.4mm\hfill}}

\let\epsilon=\varepsilon

\begin{document}

\begin{frontmatter}

\selectlanguage{francais}
\title{\textbf{Statistique sur la cyclicit\'e de Modules de Drinfeld de rang $2$}}

\selectlanguage{english}
\title{\textbf{Statistics aboute the cyclicity of a Drinfeld Modules of rank $2$}}

\author[a]{Mohamed-Saadbouh MOHAMED-AHMED}
\ead{mohamed-saadbouh.mohamed-ahmed@univ-lemans.fr}

\address[a]{D\'epartement de Math\'ematiques, Universit\'e du Maine, Avenue Olivier Messiaen,
72085 Le Mans Cedex 9.France}

\selectlanguage{francais}
\begin{abstract}
Soit $\Phi $ un $\mathbf{F}_{q}[T]$-module de Drinfeld de rang
$2$, sur un corps fini $L$, une extension de degr\'e $n$ d'un
corps fini $\mathbf{F}_{q}$. On \'etudie la cyclicit\'e de la
structure de $A$-module induite par $\Phi $ sur $L$.

{\it Pour citer cet article: Mohamed-saadbouh.Mohamed-Ahmed, C.
R. Acad. Sci. Paris,
 Ser. I ... (...).}

\vskip 0.5\baselineskip \selectlanguage{english}
\noindent{\bf Abstract} \vskip 0.5\baselineskip \noindent Let
$\Phi $ be a Drinfeld $\mathbf{F}_{q}[T]$-module of rank $2$, over
a finite field $L=\mathbf{F}_{q^{n}}$. We will study the cyclic
property of the structure $L^{\Phi }.$ We will prove that the
latter is cyclic only for trivial extensions of $\mathbf{F}_{q}$.
{\it To cite this article: Mohamed-Saadbouh Mohamed-Ahmed , C. R.
Acad. Sci. Paris, Ser. I ... (...).}
\end{abstract}

\end{frontmatter}

\vspace{-0,5 cm}

\section{Introduction}

 let $K$ a no empty global field of characteristic $p$ (
 namely a rational functions field of one indeterminate over a
 finite field ) together with a constant field, the finite field
 $\mathbf{F}_{q}$ with $p^{s}$ elements. We fix one place of $K$,
 denoted by $\infty ,$ and call $A$  the ring of regular elements
 away from the place $\infty $. Let $L$ be a commutator field of
 characteristic $p$,  $\gamma :A\rightarrow L$ be a ring
 $A$-homomorphism. The kernel of this $A$-homomorphism is denoted by $P.$
 We put $m$ =$[L,$ $A/P]$, the extension degree of $L$ over
 $A/P$, and $d=deg P$.
We denote by $L\{\tau \}$  the polynomial ring of $\tau $, namely,
the Ore polynomial ring, where $\tau $ is the Frobenius of
 $\mathbf{F}_{q}$ with the usual addition and where
 the product is given by the commutation rule: for every $%
 \lambda \in L\mathit{,}$ we have $\tau \lambda =\lambda ^{q}\tau
 $. A Drinfeld $A$-module  $\Phi :A \rightarrow
 L\{\tau \}$ is a non trivial ring homomorphism and a non trivial embedding
 of $A $ into $L\{\tau \}$ different from $\gamma $.
 This homomorphism $\Phi $, once defined, define an $A$-module structure over the $A$-field $L$, noted \ $%
 L^{\Phi }$, where the name of a Drinfeld $A$-module for a
 homomorphism $\Phi $. This structure of $A$-module depends on
 $\Phi $ and, especially, on his rank. We will make a statistic, analogue to the statistic
 for elliptic curves by Vladut in[4], about the ordinary Drinfeld $A$-modules
such that the $A$-modules $L^{\Phi }$ are cyclic, we note by
$C(d,m,q)$ the
proportion of the number (of isomorphisms of) ordinary Drinfeld \ $A$%
-modules, of rank 2 such that the A-modules structures $L^{\Phi
}$ are cyclic, this means : if we note by
$\#\{\Phi$, isomorphism, ordinary \} the number of classes of $L$%
-isomorphisms of an ordinary Drinfeld Modules of rank 2, we have :

$C(d,m,q)=\frac{\#\{\Phi, L^{\Phi } cyclic\}}{\#\{\Phi ,
isomorphism, ordinary\}}$
and we note by $C_{0}(d,m,q)$ the proportion of the number ( of
isogeny classes of) ordinary Drinfeld $A$-modules, of rank $2$ such that the $A$%
-modules $L^{\Phi }$ are cyclic, otherwise, if we note by
$\#\{\Phi $, isogeny, ordinary$\}$ the number of isogeny classes,
of ordinary Drinfeld modules of rank $2$, we have :
$C_{0}(d,m,q)=\frac{\#\{isogeny Classes of \Phi,L^{\Phi }%
{ cyclic}\}}{\#\{\Phi {, isogeny, ordinary}\}}$,
$C(d,m,q)=C_{0}(d,m,q)=1$ if and only if $m=d=1$.

This means that, to have a cyclic Drinfeld $A$-modules we must
have a trivial extension $L$, and we let think, in conjecture form, that for a big $q$ the values of $%
C(d,m,q)$ and $C_{0}(d,m,q)$ will tend to 1. \vspace{-0,5 cm}
\section{Cyclicity Statistics for the $A$-module L$^{\Phi }$}

We define $C(d,m,q)$ as been the ration of the number of
(isomorphism classes of) Drinfeld modules of rank 2 with cyclic
structure $L^{\Phi }$ to the number of $L$-isomorphisms classes of
ordinary Drinfeld modules of rank $2$, noted by $\#\{\Phi $,
isomorphism, ordinary$\}$: $C(d,m,q)=\frac {\# \{\Phi ,L^{\Phi }
cyclic \}}{\#\{\Phi, isomorphism, ordinary \}}$, as same, we
define $N(d,m,q)$ as been the ration of the number of (isogeny
classes of) Drinfeld modules of rank 2 with not cyclic structure
$L^{\Phi }$ to the number of $L$-isomorphisms isogeny of ordinary
Drinfeld modules of rank $2$, noted by $\#\{\Phi $, isogeny,
ordinary$\}$: $N(d,m,q)=\frac{\#\{\Phi ,L^{\Phi } non
cyclic\}}{\#\{\Phi, isogeny, ordinary\}}.$ We remark that : $\
0\leq C(d,m,q)$, $N(d,m,q)\leq 1$.
 J. Yu in [3], was proved that the charecteristic polynomial $P_{\Phi }$ of
a Drinfeld module $\Phi$ can be given by :
   $P_{\Phi }(X)=$ $X^{2}-cX+\mu P^{m},$ such that $\mu \in \mathbf{F}%
   _{q}^{\ast } $, and $c$ $\in A,$ where $\ \deg c\leq
   \frac{m.d}{2}$ by the Hasse-Weil analogue in this case.

Since the
no cyclicity of the structure $L^{\Phi }$ needs the fact that $\
i^{2}\mid
P_{\Phi }(1)$ and $i_{2}\mid (c-2)$, it is natural to introduce $i$ ( so $%
i_{2}$ ) in the calculus of $C(d,m,q)$ and $N(d,m,q)$. We fix the
characteristic polynomial $P_{\Phi }$, this means that we fix the
isogeny classes of $\Phi $, and we define :

\begin{definition}
We note by $n(P_{\Phi },i_{2})=\#\{\Phi : L^{\Phi
}=\frac{A}{(i_{1})}\oplus \frac{A}{(i_{2})}\}$.
\end{definition}

\begin{remark}
The number $n(P_{\Phi },i_{2})$ is equal to the number of
isomorphisms classes of $\Phi$ whole the $A$-module $L^{\Phi }\simeq \frac{A}{(i_{1})}%
\oplus \frac{A}{(i_{2})}$, in one isogeny classes, from where is
coming the correspondence between $\Phi $ and $i_{2}.$
\end{remark}

For $n(P_{\Phi }, i_{2})$ we have :

\begin{lemma}
Let $P_{\Phi }(X)=X^{2}-cX+\mu P^{m}$ be the characteristic
polynomial of an ordinary Drinfeld $A$-module $\Phi$ of rank 2,
and let $i_{2}$ be an unitary polynomial of $A$. Then if $\ $
$i_{2}\mid c-2$ we have : $n(P_{\Phi },i_{2})\geq 1$, else
$n(P_{\Phi },i_{2})=0$.
\end{lemma}

We can deduct :

\begin{corollary}
With the above notations :

$\#\{\Phi ,L^{\Phi } non cyclic\}=\sum\limits_{P_{\Phi
}}\sum\limits_{i_{2},i_{2}^{2}\mid P_{\Phi }(1)}n(P_{\Phi
},i_{2}).\#\{i_{2},i_{2}^{2}\mid P_{\Phi }(1) and i_{2}\mid
(c-2)\},$

$\#\{\Phi, L^{\Phi } cyclic\}=\sum\limits_{P_{\Phi
}}\sum\limits_{i_{2},i_{2}^{2}\mid P_{\Phi
}(1)}n(P_{\Phi},i_{2})$.$\#\{i_{2},i_{2}^{2}\nmid P_{\Phi }(1)
and i_{2}\mid (c-2)\}$,

and if we note by $n_{0}(P_{\Phi },i_{2})=$ $\#\{$isogeny classes
of $\Phi :L^{\Phi }=\frac{A}{(i_{1})}\oplus \frac{A}{(i_{2})}\}$
, we have :
$n_{0}(P_{\Phi },i_{2})=1$.
\end{corollary}

We note now by $\#\{\Phi ,$ isogeny, ordinary $\}$ the number of
isogeny classes, for an ordinary module $\Phi $, then we define :
$N_{0}(d,m,q)=\frac{\#\{isogeny classes of \Phi, L^{\Phi }%
 not cyclic \}}{\#\{\Phi ,isogeny, ordinary\}},$
the same for
$C_{0}(d,m,q)=\frac{\#\{isogeny Classes of \Phi, L^{\Phi }%
cyclic \}}{\#\{\Phi, isogeny, ordinary\}},$

We can so announce the following lemma :

\begin{lemma}
With the notations above, we have :
\begin{eqnarray*}
N_{0}(d,m,q) &=&\frac{\#\{i_{2},i_{2}^{2}\mid P_{\Phi }(1) and%
i_{2}\mid (c-2)\}}{\#\{\Phi , isogeny, ordinary\}}, \\
N(d,m,q) &=&\frac{\sum\limits_{P_{\Phi
}}\sum\limits_{i_{2},i_{2}^{2}\mid P_{\Phi }(1)}n(\Phi
,i_{2}).\#\{i_{2},i_{2}^{2}\mid P_{\Phi }(1) and i_{2}\mid
(c-2)\}}{\#\{\Phi , isomorphism, ordinary\}},
\end{eqnarray*}
$C_{0}(d,m,q)=\frac{\#\{i_{2},i_{2}^{2}\nmid P_{\Phi }(1) et
i_{2}\mid (c-2)\}}{\#\{\Phi , isogeny, ordinary\}}$,
$C(d,m,q)=\frac{\sum\limits_{P_{\Phi
}}\sum\limits_{i_{2},i_{2}^{2}\nmid P_{\Phi }(1)}n(\Phi
,i_{2}).\#\{i_{2},i_{2}^{2}\nmid P_{\Phi }(1)  and i_{2}\mid
(c-2)\}}{\#\{\Phi , isomorphism, ordinary\}} $,
and $N(d,m,q)+C(d,m,q)=1$, $N_{0}(d,m,q)+C_{0}(d,m,q)=1$.
\end{lemma}

The calculus of $\#\{\Phi ,$ isogeny, ordinary$\}$, for an ordinary $A$%
-module $\Phi $, has been calculated in [1], as been :

\begin{proposition}
Let $L=F_{q^{n}}$ and $P$ the $A$-characteristic of $L$. We put
$m=[L:A/P]$ and $d=$deg $P$ :
\begin{enumerate}
\item $m$ is odd and $d$ is odd :%
$\#\{\Phi ,isogeny, ordinary\}=(q-1)(q^{[\frac{m}{2}d]+1}-q^{[\frac{%
m-2}{2}d]+1}+1).$
\item $m.d$ even:
$\#\{\Phi , isogeny, ordinary\}=(q-1)(\frac{(q-1)}{2}q^{\frac{m}{ 2}%
d}-q^{\frac{m-2}{2}d}+1).$
\end{enumerate}
\end{proposition}

As for the number $L$-isomorphisms classes, we will need the
following result, for the proof and more details see [3] :

\begin{proposition}
Let $L$ be a finite extension of degree $\ n$ over
$\mathbf{F}_{q}$, then the number of $L$-isomorphisms classes of
a Drinfeld $A$-module of rank 2 over $L$ is $(q-1)q^{n}$ if $n$
is odd and $q^{n+1}-q^{n}+q^{2}-q$ else.
\end{proposition}

And to calculate the number of $L$-isomorphisms classes for an
ordinary Drinfeld modules, we will need to calculate the number
of $L$-isomorphisms classes for a supersingular Drinfeld modules
and subtract it from the global number of $L$-isomorphisms
classes, for this, we have by [3] :

\begin{proposition}
Let $L$ be a finite extension of $ n$ degrees over
$\mathbf{F}_{q}$, then the number of $L$-isomorphisms classes of a supersingular Drinfeld $A$%
-module of rank 2, over $L$ is $(q^{n_{2}}-1)$, where $n_{2}=$
pgcd($2,n $).
\end{proposition}

The calculus of $C(d,m,q)$ will be calculated in function of the values of $%
d $ and $m$ which are two major values to determinate $c$ because
deg $c\leq \frac{m.d}{2}$. And to calculate the number of
$L$-isomorphisms classes existing in each isogeny classes, we
need the following Definition for more information, see [4]:

\begin{definition}
Let $L$ be a finite extension of degree $\ n$ over
$\mathbf{F}_{q},$ we define $W(F)$ as been :

$W(F)=\sum\limits_{\Phi , F= Frobenius(\Phi,L) }Weigh(\Phi )$
where : $ Weigh(\Phi )=\frac{q-1}{\# Aut_{L}\Phi }.$

\end{definition}

$W(F)$ is the sum of weights( noted Weigh$(\Phi )$ ) of number of $L$%
-isomorphisms classes existing in each isogeny classes of
 the module $\Phi $ which the Frobenius is $F$. And to calculate $\#$Aut$_{L}\Phi $ we have the following lemma :

\begin{lemma}
Let $\Phi $ be an ordinary Drinfeld $A$-module of rank $2$, over
a finite field $L=F_{q^{n}} $, then : $\#$Aut$_{L}\Phi $ = $q-1$.
\end{lemma}

By the previous lemma, we can see that Weight $(\Phi )=\frac{
q-1}{\# Aut_{L}\Phi }$ $=1$, that means :

\begin{corollary}
In the case of ordinary Drinfeld modules of rank 2, $W(F)$ is the number of $%
L$-isomorphisms classes existing in each isogeny classes.
\end{corollary}

\begin{definition}
Let $D$ be an imaginary discriminant and let $l$ a polynomial for
which the square is a divisor of $D$ and let $h(\frac{D}{l^{2}})$
the number of classes of the order for which the discriminant is
$\frac{D}{l^{2}}$. We define the number of classes of Hurwitz for
an imaginary discriminant $D$, noted $H(D) $ by :
$H(D)=\sum\limits_{l}\sum\limits_{l^{2}\shortmid
D}h(\frac{D}{l^{2}}).$

\end{definition}

\begin{lemma}
If $\alpha $ is an integral element over $A$, for which
$O=A[\alpha ]$ is an $A$-order, then disc$(A[\alpha ])$ is equal
to the discriminant of the minimal polynomial of $\alpha .$
\end{lemma}

What is interesting for us is the calculus of the disc$(A[F])$ and since disc%
$(A[F])=$ disc$(P_{\Phi })$. To calculate the number of classes
$W(F)$, we have the following result, for proof see [4].

\begin{proposition}
Let $L$ be \ a finite extension of degree $n$ of a field $F_{q}$
and $F$ the Frobenius of $L$, then:

$W(F)=H(disc(A[F])).$
\end{proposition}

It remains for us to calculate $n(\Phi ,i_{2}):$

\begin{lemma}
Let $P_{\Phi }$ be the characteristic polynomial of an ordinary Drinfeld $A$%
-module of rank $2$, over a finite field $L$ such that $L^{\Phi }= \frac{A}{%
(i_{1})}\oplus \frac{A}{(i_{2})}$, and let $\Delta $ the
discriminant of the characteristic polynomial of the Frobenius
$F$, then : $n(P_{\Phi },i_{2})=H(O(\Delta /i_{2}^{2})).$
\vspace{-0,5 cm}
\end{lemma}
\section{Application}
1) $d=m=1$, in this case $L=A/P=\mathbf{F}_{q}$, the $A$-module
$L^{\Phi }$ $=A/P$ is cyclic, so $C(1,1,q)=1$. Conversely:

\begin{theorem}
Let $L=F_{q^{n}}$ and $P$ the $A$-characteristic of $L$,
$m=[L,A/P]$ and $d=$ deg$P$. Then:
$C_{0}(d,m,q)=C(d,m,q)=1\Leftrightarrow m=d=1$.

\end{theorem}

2)$m=1$ et $d=2$, In this case $n=m.d=2$, and $n_{2}=2\Rightarrow
\#\{\Phi $, isomorphism, ordinary$\}=$
$q^{3}-q-(q^{2}-1)=q^{3}-q^{2}-q+1.$
$C_{0}(2,1,q)=\frac{q(q-1)-5}{q(q-1)-2}$,

$C(2,1,q)=\frac{q^{3}-q^{2}-q+1-[\frac{q-1}{2}\sum\limits_{P_{\Phi
}}\sum\limits_{i_{2,}i_{2}^{2}\shortmid 4-4\mu P}H(O(\frac{4-4\mu P}{%
i_{2}^{2}}))+(q-1)\sum\limits_{P_{\Phi
}}\sum\limits_{i_{2},i_{2}^{2}\shortmid c^{2}-4\mu
P}H(O(\frac{c^{2}-4\mu P}{ i_{2}^{2}}))]}{q^{3}-q^{2}-q+1}$.

3)$m=2$ and $d=1$, in this case $n = m.d = 2 $, and
$n_{2}=2\Rightarrow \#\{\Phi ,$ isomorphism,
ordinary$\}=q^{3}-q-(q^{2}-1)=q^{3}-q^{2}-q+1$.
$C_{0}(1,2,q)=\frac{(q-1)q-4}{(q-1)q-2}$,
$C(1,2,q)=\frac{q^{3}-q^{2}-q+1-\sum\limits_{P_{\Phi
}}\sum\limits_{i_{2,}i_{2}^{2}\shortmid c^{2}-4\mu
P}H(O(\frac{c^{2}-4\mu P }{i_{2}^{2}}))}{q^{3}-q^{2}-q+1}.$

By the calculus of $C_{0}(d,m,q)$ and $C(d,m,q)$ for $m.d$ $\leq
2$, we have : $\lim\limits_{q\rightarrow \infty
}C_{0}(1,1,q)=\lim\limits_{q\rightarrow \infty
}C_{0}(1,2,q)=\lim\limits_{q\rightarrow \infty }C_{0}(2,1,q)=1,$
$\lim\limits_{q\rightarrow \infty
}C(1,1,q)=\lim\limits_{q\rightarrow \infty
}C(1,2,q)=\lim\limits_{q\rightarrow \infty }C(2,1,q)=1.$ By the
results above, we can give the following conjecture :

\begin{conjecture}
Let $L=F_{q^{n}}$ and $P$ the $A$-characteristic of $L$,
$m=[L,A/P]$ and $\ d= $ deg$P$. Then :

$\lim\limits_{q\rightarrow \infty
}C(d,m,q)=\lim\limits_{q\rightarrow \infty }C_{0}(d,m,q)=1$.
\end{conjecture}

\selectlanguage{english}

\end{document}